\documentclass[12pt]{amsart}

\usepackage{amssymb}

\usepackage{enumerate}

\usepackage{graphicx}

\usepackage{epstopdf}

\makeatletter
\@namedef{subjclassname@2010}{%
  \textup{2010} Mathematics Subject Classification}
\makeatother

\theoremstyle{plain}

\newtheorem*{theorem*}{Theorem}
\newtheorem*{dp*}{Dichotomy Principle (DP)}
\newtheorem*{cp*}{Containment Principle (CP)}
\newtheorem*{hfcp*}{Finite Containment Principle (FCP)}
\newtheorem*{qjcp*}{Quasi-Jordan Filling Principle (QJFP)}
\newtheorem*{hjcp*}{Jordan Filling Principle (JFP)}

\newtheorem*{DP*}{Darboux-Picard Theorem (DPT)}

\newtheorem*{fmaxx}{Finite Maximum Modulus Principle (FinMaxMP)}
\newtheorem*{maxx}{Maximum Modulus Principle (MaxMP)}

\newtheorem*{minn}{Minimum Modulus Principle (MinMP)}
\newtheorem*{prop}{Argument Containment Proposition (ACP)}

\theoremstyle{definition}

\theoremstyle{definition}
\newtheorem{ex}{Example}
\theoremstyle{remark}

\theoremstyle{remark}

\newtheorem*{remark*}{Remark}
\newtheorem*{rfp*}{Remark on the QJFP}
\numberwithin{equation}{section}

\numberwithin{equation}{section}

\newcommand{\length}{\operatorname{length}}

\newcommand{\C}{\mathbb{C}}
\newcommand{\CC}{\overline{\C}}
\newcommand{\DD}{\overline{D}}

\newcommand{\dD}{\partial D}
\newcommand{\B}{\overline{B}}
\newcommand{\z}{\overline{z}}
\renewcommand{\c}{{\operatorname{\mathsf{c}}}}

\newcommand{\arccot}{{\operatorname{arccot}}}

\newcommand{\I}{\mathcal{I}}

\newcommand{\Co}{\mathcal{C}}
\newcommand{\Fi}{\mathcal{F}}

\newcommand{\Ga}{\Gamma}
\newcommand{\De}{\Delta}

\renewcommand{\le}{\leqslant}
\renewcommand{\ge}{\geqslant}

\begin{document}


\baselineskip=17pt



\title[A topological dichotomy]{A topological dichotomy, with applications to complex analysis}

\author[I. Pinelis]{Iosif Pinelis}
\address{Department of Mathematical Sciences\\
Michigan Technological University\\
Houghton, Michigan 49931}
\email{ipinelis@mtu.edu}


\date{\today}

\begin{abstract}
Let $X$ be a compact topological space, and let $D$ be a subset of $X$. Let $Y$ be a Hausdorff topological space. Let $f$ be a continuous map of the closure of $D$ to $Y$ such that $f(D)$ is open. Let $E$ be any connected subset of the complement (to $Y$) of the image $f(\dD)$ of the boundary $\dD$ of $D$. Then $f(D)$ either contains $E$ or is contained in the complement of $E$.

Applications of this dichotomy principle are given, in particular for holomorphic maps, including maximum and minimum modulus principles, an inverse boundary correspondence, and a proof of Haagerup's inequality for the absolute power moments of linear combinations of independent Rademacher random variables. 
(A three-line proof of the main theorem of algebra is also given.)  
More generally, the dichotomy principle is naturally applicable to conformal and quasiconformal mappings. 
\end{abstract}



\maketitle

\section{Images of subsets: filling/containment dichotomy} \label{images}


Let $X$ be a compact topological space, and let $D$ be a subset of $X$, 
with the closure $\DD$ and ``boundary'' $\dD:=\DD\setminus D$; note that, if the set $D$ is open, then $\dD$ will be the boundary of $D$ in the usual sense.  
Let $Y$ be a Hausdorff topological space. 
Let $f\colon\DD\to Y$ be a continuous map such that $f(D)$ is open. 

\begin{dp*}
Let $E$ be any connected subset of the complement $f(\dD)^\c$ of $f(\dD)$ to $Y$. 
Then either $E\subseteq f(D)$ or $f(\DD)\subseteq E^\c$; that is, the image $f(D)$ of $D$ under $f$ either \emph{fills} the set $E$ or is \emph{contained} in $E^\c$.  
\end{dp*}

The proof is almost trivial. Observe first that $E\subseteq f(D)\cup f(\DD)^\c$ (since $E\subseteq f(\dD)^\c$). Next, $f(D)$ is open (by assumption) and $f(\DD)^\c$ is open as well (since $f(\DD)$ is compact). 
Therefore and because $E$ is connected, either $E\cap f(D)=E$ (that is, $E\subseteq f(D)$) or $E\cap f(\DD)^\c=E$ (that is, $f(\DD)\subseteq E^\c$). 
\qed

The DP can be rewritten in the following ``containment'' form. 
\begin{cp*}
One has $f(\DD)\subseteq\big(\bigcup_{y\in f(\DD)^\c}E_y\big)^\c$, where 
$E_y$ denotes the connected component of $y$ in $f(\dD)^\c$.  
\end{cp*}

Indeed, take any $y\in f(\DD)^\c$. Then $y\in f(\dD)^\c$ and $y\in E_y\setminus f(D)$, whence $E_y\not\subseteq f(D)$. 
So, 
by the DP, 
$f(\DD)\subseteq E_y^{\ \c}$. Thus, the DP implies the CP. 

Vice versa, suppose now that the CP holds. Let $E$ be any connected subset of $f(\dD)^\c$. 
Suppose that the first alternative, $E\subseteq f(D)$, in the DP is false. Then there exists some $y\in E\setminus f(D)$, so that $y\in f(\DD)^\c$ \big(since $y\in E\subseteq f(\dD)^\c$\big). Hence, by the CP, $f(\DD)\subseteq E_y^{\ \c}\subseteq E^\c$. 
%
\qed

The following, 
tripartite corollary of the DP may be viewed as an abstract, topological generalization of the Jordan Filling Principle for $Y=\CC$ presented in the next section. 


\begin{qjcp*}
Suppose that $D\ne\emptyset$ and let $E$ and $F$ stand for some connected subsets of $Y$. 
Then one has the following. 
\begin{enumerate}
\item[\emph{(I)}] If $f(D)\subseteq E\subseteq f(\dD)^\c$, then $f(D)=E$. 
\item[\emph{(II)}] If 
$E\subseteq f(\dD)^\c\subseteq E\cup F$,  
$f(D)\subseteq F^\c$,   
and $f(\dD)\subseteq\overline F$,  
then $f(D)=E$   
\big(moreover, it follows that 
$E\ne\emptyset$, 
$E\cap F=\emptyset$, and 
$f(D)\subseteq f(\dD)^\c$ -- that is, $f$ does not take on $D$ any of the values it takes on $\dD$\big). 
\item[\emph{(III)}] If $f(\dD)^\c=E\cup F$,  
$F\not\subseteq f(D)$, 
and $f(\dD)\subseteq\overline F$,  
then $f(D)=E$ \big(moreover, it follows that $E\ne\emptyset$, $F\ne\emptyset$, 
$E\cap F=\emptyset$, and $f$ does not take on $D$ any of the values it takes on $\dD$\big).  
\end{enumerate}
\end{qjcp*}

\begin{proof}\ \\
(I). The conditions $D\ne\emptyset$ and $f(D)\subseteq E$ imply $f(D)\not\subseteq E^\c$ and hence 
$f(\DD)\not\subseteq E^\c$. 
So, by the DP, $E\subseteq f(D)\subseteq E$, which proves part (I) of the QJFP.

\noindent (II). Assume now that the conditions of part (II) of the QJFP hold. 
Verify first the last conclusion of part (II) -- that $f(D)\subseteq f(\dD)^\c$; indeed, if that conclusion were false, then one would have
$\emptyset\ne f(D)\cap f(\dD)\subseteq f(D)\cap\overline F$, which would contradict the conditions that $f(D)$ is open and $f(D)\subseteq F^\c$. 
So, $f(D)\subseteq f(\dD)^\c\cap F^\c\subseteq(E\cup F)\cap F^\c=E\setminus F$. 
Therefore and by part (I), $E=f(D)\subseteq E\setminus F$. 
This in turn yields $E\cap F=\emptyset$. 
Also, the conditions $D\ne\emptyset$ and $E=f(D)$ imply $E\ne\emptyset$. 
This completes the proof of part (II).  

\noindent (III). Part (III) follows from part (II). Indeed, the condition $f(\dD)^\c=E\cup F$ implies $F\subseteq f(\dD)^\c$; hence, by the DP, the condition $F\not\subseteq f(D)$ yields $f(D)\subseteq F^\c$, so that all the conditions of part (II) hold. Also, the condition $F\not\subseteq f(D)$ implies $F\ne\emptyset$. 
\end{proof}

In the above proof, we deduced QJFP(II) from QJFP(I), and QJFP(III) from QJFP(II). So, one may say that QJFP(I) is the most general of the three parts of the QJFP, while QJFP(III) is the most special one.  

In the case when $f\colon\Omega\to\Omega'$ is a proper holomorphic map, where $\Omega$ and $\Omega'$ are open connected subsets of $\C^n$, Rudin \cite[Proposition~15.1.5]{rudin} shows that the ``filling'' conclusion $f(\Omega)=\Omega'$ holds  \big($f$ is said to be \emph{proper} if $f^{-1}(K)$ is compact in $\Omega$ for any compact $K\subseteq\Omega'$\big).   
Rudin \cite[Theorem~15.1.6]{rudin} also shows that, for a holomorphic map $f\colon\Omega\to\Omega'$ to be locally proper and hence open, it is enough that  
the set $f^{-1}(w)$ be compact (or, equivalently, finite) for every $w\in\Omega'$. 

More generally, the topological DP is naturally applicable to conformal and quasiconformal mappings. 

The QJFP \big(especially its parts (II) and (III)\big) will be quite useful in certain contexts, such as the proof of the JFP in the next section. 
However, at this point let us just present a simple, almost trivial illustration of how the QJFP can be applied: 

\medskip 
\emph{If $D=X\ne\emptyset$ and $Y$ is connected, then $f(D)=Y$.}

\medskip 
This follows immediately by invoking the QJFP(II) with $E=Y$ and $F=\emptyset$.

Perhaps surprisingly, the purely topological (and almost trivial) dicho\-tomy principle (DP) turns out to be 
convenient and useful in the applications to various interesting inequalities, even in the special case when the map $f$ is holomorphic.  

\section{Special cases and applications} \label{appls}

\subsection{Case $Y=\CC$} \label{special}
In this subsection, let us assume that the general conditions stated in the first paragraph of Section~\ref{images} hold. In addition, assume here that $Y=\CC:=\C\cup\{\infty\}$, the Riemann sphere, whereas $X$ may still be any compact topological space. 

However, when it is the case that $X$ equals $\CC$, $D$ is a domain (that is, an open connected set),   
and the map $f\colon\DD\to\CC$ is non-constant and holomorphic on $D$, then, by the open map theorem (cf.\ e.g.\ \cite[Theorem~5.77]{burckel} or \cite[VI.I.3]{cartan}), the condition that $f(D)$ be open will be satisfied. 
\big(Here we shall say that $f$ is \emph{holomorphic} on $D$ if for any point $z_0\in D$ there exist M\"obius transformations $M_1$ and $M_2$ of $\CC$ such that $M_1(z_0)$ is finite (that is, is in $\C$) and the function $M_2\circ f\circ M_1^{-1}$ is finite and differentiable (in the complex-variable sense) in a neighborhood of $M_1(z_0)$. Clearly, if $D\subseteq\C$ and $f(D)\subseteq\C$, then this extended notion of a holomorphic function is equivalent to the more usual one.\big)  


\begin{hfcp*}
If $f$ is finite on $\DD$,   
then $f(\DD)\subseteq
E_\infty^{\ \c}$.  
\end{hfcp*}

This follows immediately from the CP. \qed 

\begin{hjcp*}
Suppose that $f$ is finite on $\DD$ and $f(\dD)=J$, where $J$ is the image (in $\C$) of a Jordan curve. Then $f(D)=\I(J)$, where $\I(J)$ denotes the \emph{inside} of $J$, that is, the bounded connected component in $\C$ of $\C\setminus J$.  
\end{hjcp*}

This follows immediately from the QJFP(III) \big(on letting $E:=\I(J)$ and $F:=E_\infty$\big). 
\qed   

The JFP may be compared with the following result, based on the argument principle (cf.\ e.g.\ \cite[Corollary~9.16 and Exercise~9.17]{burckel}): 
\begin{DP*}
Assume that $X=\C$, $D$ is a domain, and the function $f$ is non-constant and holomorphic on $D$. 
Let $D$ be the inside of the image of a Jordan curve, and suppose that $f$ is finite on $\DD$ and 
one-to-one on $\dD$. Then $f$ is one-to-one on $\DD$, and $f(D)$ is the inside of $f(\dD)$. 
\end{DP*}
A partial extension of the DPT to $X=Y=\C\sp n$ was given by Chen \cite{chen-picard}, 
where, in addition to the injectivity of $f$ on $\DD$, 
it was proved only that $f(D)$ is \emph{a subset of} the inside (rather than \emph{exactly} the inside) of $f(\dD)$. 
  
One can see that, in contrast with the DPT, in the JFP we do not require that $D$ be a domain, or that $f$ be one-to-one on the boundary $\dD$, or that $\dD$ be the image of a Jordan curve (or any other curve), or even that the space $X$ be $\C$ or $\CC$.  On the other hand, the conclusion of the JFP is somewhat weaker than that of the DPT, in that the former is, naturally, missing the injectivity of $f$ on $\DD$. 

Of course, the QJFP is significantly more general that the JFP, even when $X=Y=\CC$ and $f$ is holomorphic on $D$. 

\begin{ex}\label{ex:zhuk}
Let $X=Y=\CC$, $D=\CC\setminus\{0,\infty\}$, $f(z)=z+1/z$ for $z\in D$, and $f(0)=f(\infty)=\infty$, so that $\dD=\{0,\infty\}$, $f(\dD)=\{\infty\}$, and $f(\dD)^\c=\C$. 
Thus, $f(\dD)$ is not the image of a Jordan curve; so, the JFP is not applicable here, and therefore the DPT is not applicable either. 
However, one can easily apply the QJFP(II) (with $E:=\C$ and $F:=\{\infty\}$), to conclude that $f(D)=\C$. Of course, in this very simple situation the same conclusion can be obtained directly, by solving a quadratic equation. \qed
\end{ex}

The following, less trivial example may be viewed as a toy model for the setting to be considered in Subsection~\ref{haag}.

\begin{ex}\label{ex:zzhuk}
Let $X=Y=\CC$ and  
$
	D=\{z\in\C\colon\Re z>0, \Im z>0\}, 
$ 
so that $\dD=\{\infty\}\cup\Ga_1\cup\Ga_2$, where $\Ga_1:=\{z\in\C\colon\Re z\ge0, \Im z=0\}$ and $\Ga_2:=\{z\in\C\colon\Re z=0, \Im z>0\}$. 
Let next 
$
	f(z)=\tfrac{2z}{z^2-1}
$ 
for $z\in\DD\setminus\{1,\infty\}$, $f(1)=\infty$, and $f(\infty)=0$.   
Then $f(\dD)=f(\Ga_1)\cup f(\Ga_2)$, $f(\Ga_1)=\{\infty\}\cup\{w\in\C\colon\Im w=0\}$, and $f(\Ga_2)=\{w\in\C\colon\Re w=0, -1\le\Im w<0\}$.  
Let now $H_+:=\{w\in\C\colon\Im w\ge0\}$, $F:=f(\Ga_2)\cup H_+$ and $E:=f(\dD)^\c\setminus F
=\C\setminus f(\Ga_2)\setminus H_+$. 
Then one can verify that the condition $f(D)\subseteq F^\c$ of the QJFP(II) holds. 
Indeed, note first that for $z\in D$ one has $f(z)=\frac{2(\z|z|^2-z)}{|z^2-1|^2}$, whence $\Im f(z)<0$, so that $f(D)\subseteq H_+^{\ \c}$. Also, for $z\in D$ one has $f(z)=\frac{2|z|^2}{z|z|^2-\z}$, whence $\Re f(z)=0$ iff $|z|=1$, in which case $|f(z)|=|\frac2{z-\z}|>1$, so that $f(D)\subseteq f(\Ga_2)^\c$. 
Thus, the condition $f(D)\subseteq F^\c$ is verified. The other conditions of the QJFP(II) are even easier to check. 
Therefore, $f(D)=E$. 

However, the JFP is not applicable here (and therefore the DPT is not applicable either), because $f(\dD)$ 
cannot be the image of a simple closed curve in $\CC$; indeed, 
$f(\Ga_1)$ is a proper closed subset 
of $f(\dD)$, and yet the set $\CC\setminus 
f(\Ga_1)$ is not connected -- cf.\ e.g.\ \cite[Exercise~4.39]{burckel}. 
One may also note the following. Suppose that $z$ traces out $\Ga_1$ from $\infty$ to $1$ to $0$, and then traces out $\Ga_2$ from $0$ to $\infty$; at that,  
$f(z)$ will first trace out the positive real semi-axis from $0$ to $\infty$, then will jump to $-\infty$ and trace out the negative real semi-axis from $-\infty$ to $0$, then will trace out the vertical segment $f(\Ga_2)$ from $0$ down to $-i$, and finally will trace out $f(\Ga_2)$ back from $-i$ to $0$. (Of course, the ``jump'' from $\infty$ to $-\infty$ is not really a jump on the Riemann sphere $\CC$.)  
Thus, $f$ is not one-to-one on $\dD$. 
This example is illustrated below.  \qed
\end{ex}

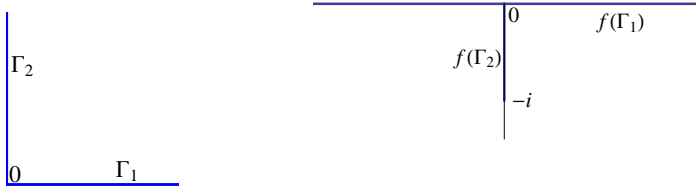
\vspace*{-1cm}
\begin{center}
\includegraphics[scale=.25]{trial
.pdf
}
\qquad\qquad
\raisebox{.6cm}{%
\includegraphics[scale=.60]{ex-right.
pdf}
}
\end{center}

\vspace*{-0cm}

\medskip
\hrule
\medskip

Consider now applications of the dichotomy principle to maximum and minimum modulus principles (again for any compact $X$). 
For any $r\in[0,\infty]$, let $B_r:=\{w\in\CC\colon|w|<r\}=\{w\in\C\colon|w|<r\}$ 
and $\B_r:=\{w\in\CC\colon|w|\le r\}$; 
one may note that the closure $\overline{B_r}$ of $B_r$ coincides with $\B_r$ unless $r=0$, 
in which latter case $\overline{B_r}=\emptyset$ and $\B_r=\{0\}$. 
Let also $M:=\sup|f|(\dD)$ and $m:=\inf|f|(\dD)$. 

\begin{fmaxx}\label{max} If $f$ is finite (on $\DD$) then \\  
$\max|f|(\DD)=\sup|f|(\dD)$. 
\end{fmaxx}

Indeed, if $M=\infty$ then the FinMaxMP is trivial. Assume now that $M<\infty$ and let 
\rule{0pt}{10pt}
$E:=\B_M^{\ \;\c}$. Then $E$ is a connected subset of $f(\dD)^\c$ and $E\not\subseteq f(D)$, since $\infty\in E\setminus f(D)$. So, by the DP, 
\rule{0pt}{10pt}
$f(\DD)\subseteq E^\c=\B_M$. \qed

More generally, the DP (with $E=\B_M^{\ \c}$) immediately yields 

\begin{maxx}\label{maxx} 
Either 
\begin{align}
& f(D)\supseteq\B_M^{\ \,\c}, \tag{$\Fi_{\max}$} \\
& \text{that is, $f$ takes on $D$ all the values that are $>M$ in modulus;} \notag \\
\intertext{or }
&f(\DD)\subseteq\B_M, \tag{$\Co_{\max}$} \\ 
&\text{that is, all the values that $f$ takes on $\DD$ are $\le M$ in modulus.} \notag  
\end{align}
Note that the ``containment'' alternative $(\Co_{\max})$ can be rewritten as \break $\max|f|(\DD)=\sup|f|(\dD)$; cf.\ the \emph{FinMaxMP}. 
\end{maxx} 

Quite similarly, the DP (with $E=B_m$) yields 

\begin{minn}\label{minn} 
Either 
\begin{align}
& f(D)\supseteq B_m, \tag{$\Fi_{\min}$} \\ 
&\text{that is, $f$ takes on $D$ all the values that are $<m$ in modulus;} \notag \\
\intertext{or }
&f(\DD)\subseteq B_m^{\ \,\c}, \tag{$\Co_{\min}$} \\ 
&\text{that is, all the values that $f$ takes on $\DD$ are $\ge m$ in modulus.} \notag 
\end{align}
Note that the ``containment'' alternative $(\Co_{\min})$ can be rewritten as \break  $\min|f|(\DD)=\inf|f|(\dD)$. 
\end{minn} 

Note also that each of the two alternatives in the MaxMP and in the MinMP (and thus in the DP) actualizes. Indeed,  take the trivial example of $f(z)=z$ for all $z\in\DD$, where $D$ is either $B_1$ or $\B_1^{\ \,\c}$. 

Various versions of the maximum and minimum modulus principles (for non-constant finite holomorphic functions on domains in $X=\C$) may be found e.g.\ in \cite{conway}. The FinMaxMP presented above corresponds to the second of the three maximum principles given in \cite[pages~124--125]{conway}. 

Our MinMP can be compared with the minimum modulus principle stated (for non-constant finite holomorphic functions on bounded domains $D$) in Exercise~1 on page~125 of \cite{conway}, which latter has the alternative $f(D)\ni0$ instead of $f(D)\supseteq B_m$; let us refer to that statement in \cite{conway} as the 0-MinMP. 
This somewhat less informative principle, 0-MinMP, is enough to obtain immediately the 
main theorem of algebra. Indeed, let $R\in(0,\infty)$ be such that $m:=\min_{|z|=R}|f(z)|>|f(0)|$, where $f$ is a given polynomial of degree $\ge1$. Then the polynomial $f$ takes on the value $0$ in $D:=B_R$, since the alternative ($|f|\ge m$ on $\DD$) cannot take place. 

One may note that (again in the case when $f$ is a non-constant holomorphic function on $D$ and $D$ is a domain) it is not hard to deduce the general MinMP from the 0-MinMP. Indeed, fix any $w_*\in B_m$. Let $g$ be a M\"obius transformation of $\CC$ leaving each of sets $B_m$, $\partial B_m$, and $\B_m^{\ \,\c}$ invariant, and such that $g(w_*)=0$. Let $h:=g\circ f$. Then $\min|h|(\dD)=m$. So, by the 0-MinMP, either $|h|\ge m$ or $h(D)\ni0$; that is, either $|f|\ge m$ or $f(D)\ni w_*$. 
However, our MinMP is more informative and directly derived. 
 
\subsection{Haagerup's inequality}\label{haag} 

Haagerup's inequalities \cite{haag} provide exact upper and lower bounds on the absolute power moments of normalized linear combinations of independent Rademacher random variables. Unfortunately, the proof given in \cite{haag} is very long and difficult. 
Nazarov and Podkorytov \cite{naz-podk} discovered a short and ingenious way to prove Haagerup's result. 
The proof in \cite{naz-podk} contains the following statement: 
\begin{prop}\label{prop:confin} The domain 
\begin{equation}\label{eq:D}
	D:=\{z\in\C\colon0<\Re z<\tfrac\pi2,\Im z>0\}\notag
\end{equation}
is mapped into the set 
\begin{equation}\label{eq:confine}
\De_p:=\{w\in\C\colon-\tfrac{\pi p}2<\arg w\le0\} \notag
\end{equation}
by the function $f$ defined by the formula
\begin{equation}\label{eq:f}
	f(z):=z^{-p}-(\pi-z)^{-p}+(\pi+z)^{-p}-(2\pi-z)^{-p}+(2\pi+z)^{-p}-\cdots, 
\end{equation}
where $1<p<2$ and the principal branch of the power function is used, so that $z^{-p}>0$ for any $z>0$; as usual, the values of the argument function $\arg$ are assumed to be in the interval $(-\pi,\pi]$; let us also assume that $\arg0=0$. 
\end{prop}
 
The authors of \cite{naz-podk} note (on page~259) that for all $z\in D$ the points $z^{-p},(\pi+z)^{-p}, 
(2\pi+z)^{-p},\dots$ are in 
$\De_p$. Then, to conclude that $f(z)\in\De_p$ for $z\in D$, they proceed to claim that the points $-(\pi-z)^{-p},-(2\pi-z)^{-p},\dots$ are also in $\De_p$; however, this claim is obviously false: if a point $z\in D$ is close (say) to $\frac\pi4$, then the arguments of the points $-(\pi-z)^{-p},-(2\pi-z)^{-p},\dots$ are close to  $-\pi\notin[-\frac{\pi p}2,0]$. 
A less stringent argument containment, with $-\pi<\arg w\le0$ in place $-\tfrac{\pi p}2<\arg w\le0$, would allow the proof in \cite{naz-podk}  to proceed. 
One may however wonder whether the more stringent ACP holds anyway. 

One may then wonder whether the ACP can be saved by simple means such as trying to prove that each of the differences $z^{-p}-(\pi-z)^{-p},(\pi+z)^{-p}-(2\pi-z)^{-p},\dots$ is in $\De_p$. However, this latter conjecture is false, even if one instead considers partial sums of these differences; e.g., the argument of the sum of the first 100 differences is 
$<-\frac{\pi p}2(1+3.5\times10^{-18})<-\frac{\pi p}2$ for $z=10^{-30}+10^{-6}i$ and $p=19/10$. 
Alternatively, one may try to consider $f(z)$ as the sum of the terms $z^{-p},-(\pi-z)^{-p}+(\pi+z)^{-p},-(2\pi-z)^{-p}+(2\pi+z)^{-p},\dots$; however, this simple trick does not work either, as already the term $-(\pi-z)^{-p}+(\pi+z)^{-p}$ is outside $\De_p$ e.g.\ when $z=\frac\pi4+10^{-2}i$ and $p=19/10$. 


Fortunately, the ACP can be rather easily proved using the topological dichotomy principle (DP). Thus, the DP method can be effective in situations where no other methods seem to work.  


\begin{proof}[Proof of the ACP]
Note that $\dD=\Ga_1\cup\dots\cup\Ga_5$, where (trying to keep up with the corresponding notation in \cite{naz-podk})
\begin{align*}
\Ga_1&:=\{0\}, \\
\Ga_2&:=\{z\in\C\colon\Re z=0,\;\Im z>0\}, \\
\Ga_3&:=\{\infty\}, \\
\Ga_4&:=\{z\in\C\colon\Re z=\tfrac\pi2,\;\Im z\ge0\}, \\
\Ga_5&:=\{z\in\C\colon0<\Re z<\tfrac\pi2,\;\Im z=0\};
\end{align*}
this is illustrated by the picture below on the left, with a portion of $D$ near $\infty$ cut off. 
Note that $|f(z)-z^{-p}|\le(\frac\pi2)^{-p/2}(1^{-p}+2^{-p}+\cdots)<\infty$ for all $z\in D$; so, 
by dominated convergence, $f$ can be extended to $\DD\setminus\{0,\infty\}$ by the same formula \eqref{eq:f}; let also $f(0):=\infty$ and $f(\infty):=0$, so that $f$ is continuous on $\DD$, and non-constant and holomorphic on $D$. 

Thus, $f=\infty$ on $\Ga_1$ and $f=0\in\De_p$ on $\Ga_3$. Let us now prove that the images of $\Ga_2$, $\Ga_4$, $\Ga_5$ under $f$ are contained in $\De_p$. \big(These images, as well as part of the boundary of the angular set $\De_p$, are shown in the picture below on the right, with a portion of $\De_p$ near $\infty$ cut off. 

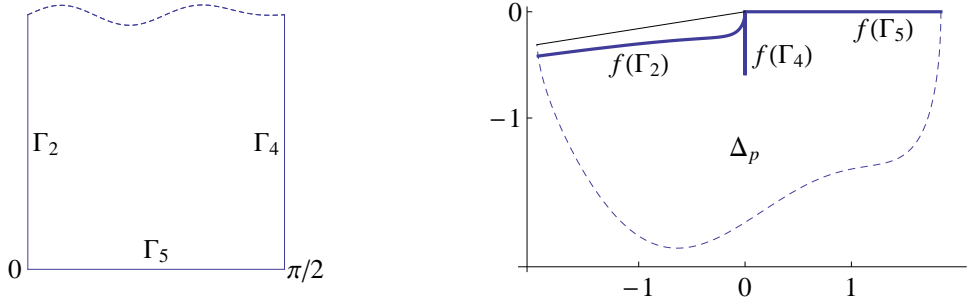
\begin{center}
\includegraphics[scale=.5]{6a.
pdf}
\qquad\qquad\quad
\raisebox{.1cm}{%
\includegraphics[scale=.75]{6.
pdf}
}
\end{center}

For $z\in\Ga_4$, one has $\Re f(z)=0$ and $\Im f(z)\le0$ (cf.\ the second displayed formula on page~262 in \cite{naz-podk}). 
Thus, $f(\Ga_4)\subseteq\De_p$. 
At this point one may note that $f(\frac\pi2)=f(\infty)=0$; so, $f(z)$ traces out the vertical segment $f(\Ga_4)$ on the imaginary axis (at least) twice as $z$ traces out the vertical ray $\Ga_4$. 
Therefore, the one-to-one condition of the Darboux-Picard theorem stated in Subsection~\ref{special} does not hold here. Yet, the dichotomy principle (DP) of Section~\ref{images} allows us to proceed and obtain the containment result.  

For $z\in\Ga_5$, one has $\Im f(z)=0$ and $\Re f(z)>0$, since $(k\pi+t)^{-p}>(k\pi+\pi-t)^{-p}$ for all $k\ge0$ and $t\in(0,\frac\pi2)$. Thus, $f(\Ga_5)\subseteq\De_p$. 

Finally, for $z\in\Ga_2$ one has $\Re f(z)<0$ and $\Im f(z)<\Re f(z)\,\tan(-\tfrac{\pi p}2)$ (cf.\ \cite[the middle of page~261]{naz-podk}). 
Thus, $f(\Ga_2)\subseteq\De_p$. 

We conclude that $f(\dD)\subseteq\De_p\cup\{\infty\}$. 
Let now $E:=\C\setminus\De_p$. Then $E$ is a connected subset of $f(\dD)^\c$. 
Moreover, every $w\in\C$ with $\Im w>0$ is in $E$. 
On the other hand, $\Im w<0$ for any $w\in f(D)$ \big(since $\Im[(k\pi+z)^{-p}]<0$ and $\Im[-(k\pi+\pi-z)^{-p}]<0$ for any $k\ge0$ and any $z\in D$\big). 
So, $E\not\subseteq f(D)$. 
Thus, by the dichotomy principle, 
$f(\DD)\subseteq E^\c=\De_p\cup\{\infty\}$. Since $f(D)\subseteq\C$, it follows that $f(D)\subseteq\De_p$. 
In fact, since $\Im w<0$ for any $w\in f(D)$, one has a little more: $f(D)$ is contained in the interior of $\De_p$. 
\end{proof}

While the ACP is enough as far as the proof of Haagerup's inequality is concerned, 
one might want to prove more. One improvement is easy. Let $A$ and $B$ denote, respectively, the sets of all points in $\C$ strictly \underline{a}bove and \underline{b}elow $f(\dD)$. Then the DP with $E=A$ (instead of $E=\C\setminus\De_p$ in the above proof of the ACP) yields $f(D)\subseteq A^\c=B\cup f(\dD)$, and the latter set is a proper subset of $\De_p$. 
Now, by the QJFP(II) of Section~\ref{images} with $E=B$ and $F=A\cup f(\Ga_4)$, one could conclude that $f(D)=B$ -- provided that one could show that $f(D)\subseteq f(\Ga_4)^\c$ (cf.\ Example~\ref{ex:zzhuk} of Subsection~\ref{special}); however, this does not appear easy to do.  

The inequality $\Im f(z)<\Re f(z)\,\tan(-\tfrac{\pi p}2)$ for $p\in(1,2)$ and $z=x+iy$ with $x\in(0,\frac\pi2)$ and $y>0$, implied by the ACP, can be rewritten as 
\begin{equation}
\sum _{k=0}^\infty 
\frac{\sin\big(p\,\arccot\frac{y}{\pi k+x}\big)}
   {\big((\pi  k+x)^2+y^2\big)^{p/2}} 
   >
   \sum _{k=0}^\infty \frac{\sin\big(p\,\arccot\frac{-y}{\pi k+\pi-x}\big)}
   {\big((\pi k+\pi-x)^2+y^2\big)^{p/2}} 
   \tag{*}	
\end{equation}
\big(with the values of the function $\arccot$ in the interval $(0,\pi)$\big).    
To appreciate the usefulness of the topological dichotomy principle, one may try to prove (*) by other methods, say by methods of the calculus of functions of real variables.


%
%
%
\end{document}